\documentclass{amsart}

\usepackage{xcolor,graphicx}
\usepackage{hyperref}
\usepackage{tikz}
\usepackage{amsmath, amsfonts, latexsym, array, amssymb, amscd }
\usepackage{enumitem}
\usetikzlibrary{shapes,arrows,positioning,decorations.pathreplacing}

\usepackage{ulem}

\numberwithin{equation}{section}

\theoremstyle{theorem}
\newtheorem{theorem}{Theorem}[section]
\newtheorem{proposition}[theorem]{Proposition}

\address{Z.~Li, Suzhou High School of Jiangsu Province, No 699 Renmin Road, Suzhou 215007, China, \emph{E-mail address:} \tt{2209313396@qq.com}}

\address{M.~Lu, Suzhou High School of Jiangsu Province, No 699 Renmin Road, Suzhou 215007, China, \emph{E-mail address:} \tt{2528326549@qq.com}}

\address{T.~Wang, School of Mathematical Sciences, Shanghai Jiaotong University, No 800 Dongchuan Road, Shanghai 200240, P.R.China, \emph{E-mail address:} \tt{a356789xe@sjtu.edu.cn}}

\begin{document}

\title[Entropy Monotonicity]{An inequality in real Milnor-Thurston monotonicity problem}

\author{Ziyu Li}
\author{Minyu Lu}
\author{Tianyu Wang}

\date{\today}
\thanks{TW was supported by NSFC grants 24Z033004105}

\maketitle

\begin{abstract}
In late 1990's, Tsujii proved monotonicity of topological entropy of real quadratic family $f_c(x)=x^2+c$ on parameter $c$ by proving an inequality concerning orbital information of the critical point. In this paper, we consider a weak analog of such inequality for the general family $f_{c,r}(x)=|x|^r+c$ with rational $r>1$, by following an algebraic approach.
\end{abstract}

\section{Preliminaries}

Given a family of $c$-parametrized continuous maps $\{f_c\}$ on $\mathbb{R}$ with fixed amount of turning points, much attention has been driven to better understand the relation between parameter $c$ and certain orbital complexity of $f_c$. Such complexity, being widely known as \emph{topological entropy} nowadays, was shown in \cite{MS77,MS80} to have a more intuitive interpretation of being the exponential growth rate of the number of \emph{laps} (i.e. interval of monotonicity) for $f_c^n$, denoted by $l(f_c^n)$. That is,
$$
h_{\text{top}}(f_c)=\lim_{n\to \infty}\frac{1}{n}\log l(f_c^n),
$$
where $h_{\text{top}}(f_c)$ refers to topological entropy of $f_c$. The study on whether $h_{\text{top}}(f_c)$ varies monotonically with $c$ traces back to 1970's and gives rise to a collection of celebrated works for a variety of systems, among which one of the most well-known model is the real quadratic family of unimodal maps. To be precise, by specializing to the Douady-Hubbard normal form 
$$
f_c(x):=x^2+c,
$$
its topological entropy is a monotone decreasing function of parameter $c$. Different proofs of this result (in the form of quadratic logistic family) are given independently using ideas from holomorphic dynamics; see Milnor-Thurston \cite{MT88}, Douady-Hubbard \cite{DH84,DH85-1,Dou95}, and Sullivan \cite[Theorem VI.4.2]{WS93}. In the mean time, Tsujii provided a completely different proof in \cite{TSU00}. The spirit of his argument lies in the verification of the following inequality concerning orbital information of the critical point
\begin{equation} \label{eqoriginalinequality}
    \frac{(f_c^{n-1})'(t)|_{t=c}}{D_c(f_c^n(0))}>0,
\end{equation}
where $(c,n)$ refers to the pairs in $\mathbb{R}\times \mathbb{N}$ such that $0$ is periodic under $f_c$ with primitive period being $n$. We call such $(c,n)$ a \emph{periodic pair}. In this paper, we consider a weak analog of \eqref{eqoriginalinequality} over a wider range of unimodal real maps $\{f_{c,r}\}$ with parameter $c\in \mathbb{R}$, defined by
$$
f_{c,r}(x):=|x|^r+c,
$$
where $r>1$ is a rational number. Our main goal is to show that for all periodic pair $(c,n)$, the following transversality condition holds true
\begin{equation} \label{eqourinequality}
    \frac{{(f^{n-1}_{c,r}})'(t)|_{t=c}}{D_c(f_{c,r}^n(0))}\neq 0.
\end{equation}
It follows immediately from the definition of primitive period that $(f_{c,r}^{n-1})'(t)|_{t=c}\neq 0$. Therefore, \eqref{eqourinequality} is equivalent to
\begin{equation} \label{eqourinequality1}
    D_c(f_{c,r}^n(0)) \neq 0.
\end{equation} 
When $r$ is fixed, we will omit $r$ in the notation and abbreviate $f_{c,r}$ as $f_c$. The main theorem is then stated as follows
\begin{theorem} \label{thmmain}
    Given any fixed rational $r>1$, for each periodic pair $(c,n)$, i.e. $n\geq 1$, $f_c^n(0)=0$, and $f_c^j(0)\neq 0$ for all $j\in [1,n-1]$, \eqref{eqourinequality1} holds true. 
\end{theorem}
Unlike the proof for \eqref{eqoriginalinequality} from Tsujii using quadratic differentials (see also \cite{Mil00} for some insights), our proof is based on an algebraic observation made from \cite[Lemma 1, page 333]{DH85}, which in our situation states that
\begin{proposition} \label{propalgebraicinteger}
    Given any fixed rational $r>1$, if $c$ is an algebraic integer, then \eqref{eqourinequality1} holds true as long as $n$ satisfies $f_c^j(0)\neq 0$ for all $j\in [1,n-1]$. 
\end{proposition}
As a consequence of Proposition \ref{propalgebraicinteger}, we are able to conclude Theorem \ref{thmmain} as long as we have the key ingredient stated as follows
\begin{proposition} \label{propkey}
    Fixing any $r>1$, for each periodic pair $(c,n)$, $c$ is an algebraic integer.
\end{proposition}

\subsection*{Acknowledgments} This project is part of Tencent Aspiring Explorers In Science Program. We would like to thank Weixiao Shen for suggesting this problem, and for proposing numerous helpful advices on improving the manuscript. We also acknowledge the hospitality of Southern University of Science and Technology, where some of this work was completed. 

\section{Proof of Proposition \ref{propalgebraicinteger} \& \ref{propkey}}

Throughout the rest of this section we fix $r=p/q$ with $p>q$ being co-prime positive integers, and write $f_c:=f_{c,r}$ as no ambiguity shall raise. We also assume WLOG that $c<0$ as otherwise there is nothing to prove.

\subsection{Proof of Proposition \ref{propalgebraicinteger}}
Fix any $c$ and $n$ be as in the statement of the proposition. Let $\xi_j:=f_c^j(0)$ for all $j\in [1,n]$, and $s_j:=\text{sgn}(\xi_j)\in \{-1,1\}$ be the sign of $\xi_j$ for all $j\in [1,n-1]$. Such $s_j$ are all well-defined due to our assumption on $c$ and $n$. It follows from a straightforward induction process that each $\xi_j$ is an algebraic integer. Also notice that under the assumption of Proposition \ref{propalgebraicinteger}, each $s_j$, being a function of $c$, is locally constant near $c$.

\emph{Claim: for any $j\in [1,n]$, $q^{j-1}D_c(f_c^j(0))=q^{j-1}D_c(\xi_j)$ is an algebraic integer.}

We first show how Proposition \ref{propalgebraicinteger} follows from the claim. It suffices to prove the case of $n\geq 2$. Assuming the claim, we know 
$$
\begin{aligned}
    D_c(\xi_n)
    &=D_c((s_{n-1}\xi_{n-1})^r)+1
=r(s_{n-1}\xi_{n-1})^{r-1}D_c(s_{n-1}\xi_{n-1})+1 \\\
&=\frac{p(s_{n-1}\xi_{n-1})^{r-1}(q^{n-2}D_c(s_{n-1}\xi_{n-1}))}{q^{n-1}}+1=\frac{pA_{n-1}}{q^{n-1}}+1,
\end{aligned}
$$
where $A_n:=(s_{n-1}\xi_{n-1})^{r-1}(q^{n-2}D_c(s_{n-1}\xi_{n-1}))$. Notice that both $(s_{n-1}\xi_{n-1})^{r-1}$ and $q^{n-2}D_c(s_{n-1}\xi_{n-1})$ are algebraic integers, provided $c$ being an algebraic integer and the claim with $j=n-1$ respectively. Therefore, $A_{n-1}$ is also an algebraic integer. If we assume $D_c(\xi_n)= 0$, it turns out that $A_{n-1}=-\frac{q^{n-1}}{p}$, which is not an algebraic integer. This contradicts the observation we just made, thus concludes the proposition provided the claim.

It remains to prove the claim is true, and we prove by induction. The case for $j=1$ is obvious as $D_c(\xi_1)=1$. Assume the claim holds for all $j\leq k$, with $k\in [1,n-1]$. Taking $j=k+1$, we have
$$
\begin{aligned}
q^kD_c(\xi_{k+1})=q^k(D_c((s_k\xi_k)^{r})+1)=p(s_k\xi_k)^{r-1}(q^{k-1}D_c(s_{k-1}\xi_{k-1}))+q^k,
\end{aligned}
$$
which implies that $q^kD_c(\xi_{k+1})$ is also an algebraic integer as both $(s_k\xi_k)^{r-1}$ and $q^{k-1}D_c(s_{k-1}\xi_{k-1})$ are algebraic integers, where the latter follows from our induction hypothesis. This concludes the proof of our claim, thus the proposition.

\subsection{Proof of Proposition \ref{propkey}}
Let $(c,n)$ be a periodic pair for $f_c=f_{c,r}$. As in the proof of the above proposition, we get rid of absolute value in the presentation of $f_c$ by introducing a (finite) sequence of symbols. For each $j\in [1,n]$, we let $b_j$ be such that
$$
f_c^j(0)=b_jc.
$$
For instance, $b_1=1$, $b_2=-(-c)^{r-1}+1$ as $f_c^2(0)=|c|^r+c=(-(-c)^{r-1}+1)c$. For each $j\in [1,n-1]$, let
$$
s_j:=\text{sgn}(b_{j}c)\in \{-1,1\}
$$
be the sign of $b_{j}c$. Then
$$
b_{j+1}c=f_c^{j+1}(0)=|f_c^j(0)|^r+c
=(s_{j}b_jc)^r+c=(-(-s_{j}b_j)^r(-c)^{r-1}+1)c
$$
which implies that
\begin{equation} \label{eqbiterate1}
    b_{j+1}=-(-s_{j}b_j)^r(-c)^{r-1}+1.
\end{equation}
Writing $s:=-(-c)^{r-1}$, \eqref{eqbiterate1} turns into
\begin{equation} \label{eqbiterate}
    b_{j+1}=(-s_{j}b_j)^rs+1.
\end{equation}
Proving Proposition \ref{propkey} is equivalent to show that $s$ is an algebraic integer. To prove it, start with noticing that $0=f_c^n(0)=b_nc$, which implies that $b_n=0$. By applying \eqref{eqbiterate} we have
$$
0=b_n=(-s_{n-1}b_{n-1})^rs+1,
$$
which gives
$$
(-s_{n-1}b_{n-1})^ps^q-(-1)^q=0.
$$
Another application of \eqref{eqbiterate} shows that
\begin{equation} \label{eqstep1}
    (-s_{n-1}((-s_{n-2}b_{n-2})^{p/q}s+1))^ps^q-(-1)^q=0.
\end{equation}
If we write $b'_{n-2}:=(-s_{n-2}b_{n-2})^{1/q}=|b_{n-2}|^{1/q}$, \eqref{eqstep1} turns into
\begin{equation} \label{eqstep1'}
    -s_{n-1}^p\sum_{i=0}^p\binom{p}{i}(b'_{n-2})^{pi}s^{q+i}-(-1)^q=0.
\end{equation}
Observe that LHS of \eqref{eqstep1'} is a multi-variable polynomial with variables being $s$ and $b'_{n-2}$. Meanwhile, this polynomial is \emph{monic} in the sense of having the coefficient for the term with highest degree (in both $s$ and $b'_{n-2}$) be equal to $\pm 1$. The main strategy to proceed from \eqref{eqstep1} is that by applying \eqref{eqbiterate} once again and conducting an appropriate manipulation to get rid of the fraction exponent $p/q$, we end with a polynomial in $s$ and $b'_{n-3}$ which is also monic in the above sense. In particular, the index of $b'$-term decreases by $1$. After repeating such process for $n-3$ times, we are left with a monic polynomial with a single variable $s$, as now $b_2$ is just $s+1$.

Now let us proceed to the construction of such manipulation. Write $p^2=k_1q+k_2$, where $k_1,k_2\in \mathbb{N}$ and $k_2\in [0,q-1]$. Also write $\{0,1,\cdots,(k_1+1)q-1\}=\biguplus_{j=0}^{q-1}I_j$, where $I_j:=\{kq+j:0\leq k \leq k_1\}$ for each $i$. By multiplying $b_{n-2}'^{q-k_2-1}$ on both sides of \eqref{eqstep1'}, we have
\begin{equation} \label{eqsumsimplified}
    \sum_{j=0}^{q-1}a_{n-2,j}b_{n-2}'^j=0
\end{equation}
where for each $j$ we write
$$
a_{n-2,j}=\sum_{i:pi+q-k_2-1\in I_j}-s_{n-1}^p\binom{p}{i}(b'_{n-2})^{pi+q-k_2-1-j}s^{q+i}.
$$
Notice that $\{a_{n-2,j}\}_{j=0}^{q-1}$ is a sequence of multi-variable polynomials in $s$ and $b_{n-2}'^q$. It is clear by plugging in $i=p$ and recalling the choice on $k_2$ that $a_{n-2,q-1}$ is monic, and has highest degree on both $b_{n-2}'^q$ and $s$ among all $\{a_{n-2,j}\}_{j=0}^{q-1}$. The main idea is to gradually `cancel' terms in \eqref{eqsumsimplified} without changing the monic behavior of the leading term. Precisely, our canceling process aims at using \eqref{eqsumsimplified} to produce
\begin{equation} \label{eqstep2}
    \sum_{j=0}^{q-2}a_{n-2,j}^{(2)}b_{n-2}'^j=0,
\end{equation}
where $\{a_{n-2,j}^{(2)}\}_{j=0}^{q-2}$ is a sequence of multi-variable polynomial in $s$ and $b_{n-2}'^q$, such that $a_{n-2,q-2}^{(2)}$ is monic, and has the highest degree on both variables.
\begin{proof}[Proof of \eqref{eqstep2}]
    We may rewrite \eqref{eqsumsimplified} as
    \begin{equation} \label{eqsumsimplified2}
        a_{n-2,q-1}b'^{q-1}_{n-2}=-\sum_{j=0}^{q-2}a_{n-2,j}b_{n-2}'^j.
    \end{equation}
    By multiplying both sides of \eqref{eqsumsimplified} by $a_{n-2,q-1}b'_{n-2}$, we have
    \begin{equation} \label{eqpf1}
        a_{n-2,q-1}^2b_{n-2}'^q+\sum_{j=0}^{q-2}a_{n-2,j}a_{n-2,q-1}b_{n-2}'^{j+1}=0.
    \end{equation}
    Plugging \eqref{eqsumsimplified2} into \eqref{eqpf1}, we have
    \begin{equation} \label{eqpf2}
        a_{n-2,q-1}^2b_{n-2}'^q+\sum_{j=0}^{q-2}\hat{a}_{n-2,j}^{(2)}b_{n-2}'^{j}=0,
    \end{equation}
    where  $\hat{a}_{n-2,j}^{(2)}:=a_{n-2,j-1}a_{n-2,q-1}-a_{n-2,j}a_{n-2,q-2}$. Apparently, we have
    \begin{itemize}
        \item Each element in $\{\hat{a}_{n-2,j}^{(2)}\}_{j=0}^{q-2}\cup \{a_{n-2,q-1}^2\}$ is a multi-variable polynomial in $s$ and $b_{n-2}'^q$.
        \item Among them, $a_{n-2,q-1}^2$ is the only term with the highest degree in both variables.
    \end{itemize}
    Multiplying both sides of \eqref{eqpf2} by $a_{n-2,q-1}b'_{n-2}$ and plugging in \eqref{eqsumsimplified2} once again, we obtain
    \begin{equation} \label{eqpf3}
        a_{n-2,q-1}^3b_{n-2}'^{q+1}+\sum_{j=0}^{q-2}\hat{a}_{n-2,j}^{(3)}b_{n-2}'^{j}=0,
    \end{equation}
    where $\hat{a}_{n-2,j}^{(3)}:=a_{n-2,q-1}\hat{a}_{n-2,j-1}^{(2)}-a_{n-2,j}\hat{a}_{n-2,q-2}^{(2)}$ for all $j\in [1,q-2]$ and $\hat{a}_{n-2,0}^{(3)}:=-a_{n-2,0}\hat{a}_{n-2,q-2}^{(2)}$. By repeating the above process for an additional $q-3$ times and a straightforward induction argument, we know 
    \begin{equation} \label{eqpf4}
        a_{n-2,q-1}^qb_{n-2}'^{2q-2}+\sum_{j=0}^{q-2}\hat{a}_{n-2,j}^{(q)}b_{n-2}'^{j}=0,
    \end{equation}
    where 
    \begin{itemize}
        \item Each element in $\{\hat{a}_{n-2,j}^{(q)}\}_{j=0}^{q-2}\cup \{a_{n-2,q-1}^q\}$ is a multi-variable polynomial in $s$ and $b_{n-2}'^q$.
        \item Among them, $a_{n-2,q-1}^q$ is the only term with the highest degree in both variables.
    \end{itemize}
    It then follows from \eqref{eqpf4} that
    \begin{equation} \label{eqpf5}
        (a_{n-2,q-1}^qb_{n-2}'^{q}+\hat{a}_{n-2,q-2}^{(q)})b_{n-2}'^{q-2}+\sum_{j=0}^{q-3}\hat{a}_{n-2,j}^{(q)}b_{n-2}'^{j}=0,
    \end{equation}
    where $a_{n-2,q-1}^qb_{n-2}'^{q}+\hat{a}_{n-2,q-2}^{(q)}$ is monic and leads in degree for both variables. Therefore, writing $a_{n-2,q-2}^{(2)}:=a_{n-2,q-1}^qb_{n-2}'^{q}+\hat{a}_{n-2,q-2}^{(q)}$ and $a_{n-2,j}^{(2)}:=\hat{a}_{n-2,j}^{(q)}$ for all other $j$, we have managed to obtain \eqref{eqstep2} satisfying our wanted properties. 
\end{proof}

The process now continues by repeating the canceling process displayed in the above proof for an additional $q-3$ times. This will provides us with the following
\begin{equation} \label{eqstep3}
    a_{n-2,1}^{(q-1)}b_{n-2}'+a_{n-2,0}^{(q-1)}=0,
\end{equation}
which implies that
\begin{equation} \label{eqstep3'}
    (a_{n-2,1}^{(q-1)})^qb_{n-2}'^q-(-a_{n-2,0}^{(q-1)})^q=0,
\end{equation}
where both $a_{n-2,1}^{(q-1)}$ and $a_{n-2,0}^{(q-1)}$ are multi-variable polynomials in $s$ and $b_{n-2}'^q$, with $a_{n-2,1}^{(q-1)}$ being monic and dominating in degree on both variables. Plugging $b_{n-2}'=(-s_{n-1}b_{n-2})^{1/q}$ into \eqref{eqstep3'}, we have
\begin{equation} \label{eqend1}
    b_{n-2}^{m_1}s^{n_1}+\sum_{i\in [0,m_1-1],j\in [0,n_1-1]}c^{1}_{i,j}b_{n-2}^is^j=0,
\end{equation}
where $m_1$ and $n_1$ are two positive integers, and $c^{1}_{i,j}\in \mathbb{Z}$ for each $(i,j)\in [0,m_1-1]\times [0,n_1-1]$. Applying \eqref{eqbiterate} for $j=n-3$ to \eqref{eqend1} and repeating the above argument  leading towards \eqref{eqend1}, we end up with
\begin{equation} \label{eqend2}
    b_{n-3}^{m_2}s^{n_2}+\sum_{i\in [0,m_2-1],j\in [0,n_2-1]}c^{2}_{i,j}b_{n-3}^is^j=0,
\end{equation}
where $m_2$ and $n_2$ are two positive integers, and $c^{2}_{i,j}\in \mathbb{Z}$ for each $(i,j)\in [0,m_2-1]\times [0,n_2-1]$. Such process can go on for an additional $n-5$ times, and the last equation we have is
\begin{equation} \label{eqend3}
    b_2^{m_{n-3}}s^{n_{n-3}}+\sum_{i\in [0,{m_{n-3}-1],j\in [0,{n_{n-3}-1]}}}c^{n-3}_{i,j}b_{n-3}^is^j=0,
\end{equation}
where $m_{n-3}$ and $n_{n-3}$ are two positive integers, and $c^{n-3}_{i,j}\in \mathbb{Z}$ for each $(i,j)\in [0,m_{n-3}-1]\times [0,n_{n-3}-1]$. Plugging in $b_2=s+1$, it follows immediately that $s$ is an algebraic integer, concluding the proof of Proposition \ref{propkey}.

\bibliographystyle{amsalpha}
\bibliography{Tencent}
\end{document}